\documentclass[a4paper,11pt]{amsart}
\usepackage[utf8]{inputenc}
\usepackage{amsmath, amssymb, amsthm}
\usepackage{geometry}
\usepackage{enumitem}
\usepackage{booktabs}
\usepackage{hyperref}

\newtheorem{theorem}{Theorem}[section]
\newtheorem{lemma}[theorem]{Lemma}
\newtheorem{proposition}[theorem]{Proposition}

\theoremstyle{definition}

\newcommand{\Q}{\mathbb{Q}}

\title[Sums of powers of two $k$-Pell numbers]{On the exponential Diophantine equation
$\left(P_{n}^{(k)}\right)^{x}+\left(P_{r}^{(k)}\right)^{x}=P_{m}^{(k)}$}

\author[F. Faye]{Fallou Faye}
\address{F. Faye, UFR SAT, Universit\'e Alioune Diop, Bambey, 30, S\'en\'egal}
\email{falloufaye1993@gmail.com}

\author[B. Faye]{Bernadette Faye}
\address{B. Faye, UFR SAT, Universit\'e Alioune Diop, Bambey, 30, S\'en\'egal}
\email{bernadette.faye@uadb.edu.sn}

\subjclass[2020]{11B39, 11D61, 11J86}
\keywords{$k$-generalized Pell numbers, exponential Diophantine equations, linear forms in
logarithms, reduction method}

\begin{document}

\begin{abstract}
The Pell numbers satisfy $P_n^2+P_{n+1}^2=P_{2n+1}$, and Bravo, Luca and Pomeo recently found all
the cases in which the sum of the $x$-th powers of two consecutive $k$-generalized Pell numbers is
again a $k$-generalized Pell number. We ask what happens when the two indices are arbitrary, that
is, we solve $\big(P_n^{(k)}\big)^x+\big(P_r^{(k)}\big)^x=P_m^{(k)}$ in positive integers $k\ge2$,
$n\ge r$, $m$ and $x$. The answer is that nothing new appears: apart from $1^x+1^x=2$, every
solution has consecutive indices. The proof combines linear forms in logarithms with an auxiliary
Fibonacci equation which does not depend on $k$, and, for large $k$, with the fact that the
dominant root of the sequence is very close to $\varphi^2$. For the reduction step,
the method of Dujella and Peth\H{o} fails for large indices and has to be replaced by a criterion
of Legendre type.
\end{abstract}

\maketitle
\section{Introduction}

Among the many identities satisfied by the Pell numbers $0,1,2,5,12,29,70,169,\ldots$, one of the
most pleasant is
\begin{equation}\label{eq:pellid}
P_n^2+P_{n+1}^2=P_{2n+1}
\end{equation}
meaning that the sum of the squares of two consecutive Pell numbers is again a Pell number, as in
$2^2+5^2=29$ or $5^2+12^2=169$. The Fibonacci numbers have the same property,
$$F_n^2+F_{n+1}^2=F_{2n+1},$$and it is natural to ask what happens with higher powers. Marques and
Togb\'e \cite{MT2010} showed that if $F_n^x+F_{n+1}^x$ is a Fibonacci number for all sufficiently
large $n$, then $x\in\{1,2\}$; Luca and Oyono \cite{LO2011} then proved that $F_n^x+F_{n+1}^x=F_m$
has no solution at all with $n\ge2$ and $x\ge3$, and G\'omez and Luca \cite{GL2014} extended this
to the $k$-generalized Fibonacci numbers. For the Pell numbers the question was settled in
\cite{RFLT2019}, by Rihane, B.~Faye, Luca and Togb\'e: $P_n^x+P_{n+1}^x=P_m$ has only the solutions
given by \eqref{eq:pellid}. For the $k$-generalized Pell numbers, Hasanalizade \cite{Has2024}
treated the squares, and Bravo, Luca and Pomeo \cite{BLP2025} recently solved the equation
$$\big(P_n^{(k)}\big)^x+\big(P_{n+1}^{(k)}\big)^x=P_m^{(k)}$$ completely.

In a companion paper \cite{FF2026} we study the equation
$$\big(P_{n+1}^{(k)}\big)^x-\big(P_{n-1}^{(k)}\big)^x=2P_m^{(k)}$$ which comes from two other Pell
identities, $$P_{n+1}-P_{n-1}=2P_n \quad \text{and}\quad  P_{n+1}^2-P_{n-1}^2=2P_{2n}.$$ While working on it, we asked
ourselves what is really special about the consecutive indices in \eqref{eq:pellid}. Could a sum
such as $P_7^x+P_3^x$, with indices far from each other, be a $k$-Pell number? For the Fibonacci numbers, the
analogous question with an arbitrary gap was considered recently by Alvarenga, Chaves, Ramos, Silva
and Sosa \cite{ACRSS2025}. In this paper we answer it for the $k$-generalized Pell numbers. Let
$k\ge2$ and let $P^{(k)}=(P_n^{(k)})_n$ be the $k$-generalized Pell sequence.. We study
\begin{equation}\label{eq:main}
\left(P_{n}^{(k)}\right)^{x}+\left(P_{r}^{(k)}\right)^{x}=P_{m}^{(k)}
\end{equation}
in positive integers $k\ge2$, $n\ge r$, $m$ and $x$, with no condition on the gap $n-r$.

We began with a computer search, for $2\le k\le80$, $n\le45$ and $x\le25$. It found nothing new.
Apart from the trivial $1^x+1^x=2$, every solution it produced had $r=n-1$ and was already in the
list of \cite{BLP2025}; the curious solution $2^5+1^5=33=P_5^{(3)}$, for example, appears for
$k=3$ only, and it is the only non-trivial solution with an exponent larger than $2$. The main result of the
paper is that the search did not miss anything.

\begin{theorem}\label{thm:main}
All the solutions of \eqref{eq:main} in positive integers $k\ge2$, $n\ge r\ge1$, $m$ and $x$ are
\begin{enumerate}[label=\rm(\alph*)]
\item $(n,r,m,x)=(1,1,2,x)$ for every $k\ge2$ and $x\ge1$, namely $1^x+1^x=2$;
\item $(n,r,m,x)=(2,1,3,2)$ for every $k\ge2$, namely $2^2+1^2=5$;
\item $(n,r,m,x)=(n,n-1,2n-1,2)$ for $k=2$ and every $n\ge2$, namely \eqref{eq:pellid};
\item $(n,r,m,x)=(2,1,5,5)$ for $k=3$, namely $2^5+1^5=33=P_5^{(3)}$.
\end{enumerate}
\end{theorem}

So the answer to our question is negative. Part of the explanation is a fact observed by K{\i}l{\i}\c{c} in 
\cite{Kilic2008}. He show that the first $k+1$ non-zero terms of $P^{(k)}$ are the Fibonacci numbers of odd
index, $P_j^{(k)}=F_{2j-1}$ for $1\le j\le k+1$. As long as all the indices in \eqref{eq:main} are
at most $k+1$, the equation does not see $k$ at all, and we are led to the following statement,
which we prove first.

\begin{theorem}\label{thm:fib}
The equation
\begin{equation}\label{eq:fib}
F_{2n-1}^x+F_{2r-1}^x=F_{2m-1}
\end{equation}
has no solution in integers $n\ge3$, $1\le r\le n$, $m\ge1$ and $x\ge2$.
\end{theorem}

The general strategy is the one of \cite{BLP2025} and of our companion paper, but the gap $d=n-r$ is
a new unknown, and this changes the argument in several places. When $d$ is large, $P_r^x$ is
negligible compared with $P_n^x$, and a linear form in the logarithms of $g$, $\alpha$ and $P_n$ is
very small. When $d$ is small, $P_r^x$ is not negligible, and we have to keep the factor
$1+\alpha^{-dx}$ in a second form; this form can vanish, and for $k=2$ it does vanish, precisely
because of \eqref{eq:pellid}. The case $r=n$ gives the equation $2P_n^x=P_m$; to see that the
corresponding form does not vanish we needed a lower bound for the discriminant of the
characteristic polynomial of $P^{(k)}$ ( see our Lemma~\ref{lem:disc}).

There was also an obstacle we did not expect. For large $n$, $\log P_n$ is extremely close to
$\log g+n\log\alpha$ , so that the
inhomogeneous part of the linear form is almost a combination of $1$ and of the irrational number
used in the reduction. The method of Dujella and Peth\H{o} then fails. We then replaced it, for these $n$, by a criterion
of Legendre type in which the small perturbation is controlled.

For large $k$ the dominant root $\alpha$ of $P^{(k)}$ is very close to $\varphi^2$, and all the
linear forms can be replaced by forms in $\Q(\sqrt5)$ whose coefficients do not depend on $k$. This
gives $k<1.57\times10^{45}$, and a few rounds of reduction bring this down to $k\le604$.  Table~\ref{tab:roadmap} summarizes
the structure of the proof.

\begin{table}[ht]
\centering
\caption{Structure of the proof.}\label{tab:roadmap}
\begin{tabular}{@{}lll@{}}
\toprule
Range & Tool & Result \\
\midrule
$x=1$, $n\le2$, $(k,x,d)=(2,2,1)$ & elementary, \cite{BT2024,NRT2021} & Section~\ref{sec:elem}\\
$n\le k+1$, $m\le k+1$ & Theorem~\ref{thm:fib} (forms in $\Q(\sqrt5)$) & no solution\\
all $k$ & forms in $\Q(\alpha)$, Matveev & $x<X(k)$, Prop.~\ref{prop:X}\\
$k\ge K_{\rm abs}$ & forms in $\Q(\sqrt5)$ & none, Prop.~\ref{prop:Kabs}\\
$K_0<k<K_{\rm abs}$ & reduction in $\Q(\sqrt5)$ & none, Prop.~\ref{prop:rounds}\\
$2\le k\le K_0$ & reduction in $\Q(\alpha)$, search & Prop.~\ref{prop:smallk}\\
\bottomrule
\end{tabular}
\end{table}

Here $K_{\rm abs}=1.5694\times10^{45}$ and $K_0=604$. All the computations were done in Python, the
details are in Appendix~\ref{app:comp}.
\section{Preliminaries}\label{sec:prelim}

\subsection{The generalized Pell sequence}

Let $k \geq 2$ be an integer. The $k$-generalized Pell sequence $P^{(k)}:=(P^{(k)}_n)_{n\geq
2-k}$ is defined by
$$P^{(k)}_n=2P^{(k)}_{n-1}+P^{(k)}_{n-2}+\cdots + P^{(k)}_{n-k}\quad\hbox{for all $n\geq 2$},$$
with the initial conditions $P^{(k)}_{-(k-2)}=P^{(k)}_{-(k-3)}=\cdots=P^{(k)}_{-1}=P^{(k)}_0=0$
and $P^{(k)}_1=1$. For $k=2$ we obtain the Pell sequence $0,1,2,5,12,29,70,\ldots$. The terms of
this sequence are called $k$-generalized Pell numbers. When $k$ is fixed we simply write
$P_n:=P_n^{(k)}$.

K{\i}l{\i}\c{c} \cite{Kilic2008} showed that the first $k+1$ non-zero terms of $P^{(k)}$ are
Fibonacci numbers with odd index, namely
\begin{equation}\label{eq:kilic}
P_n^{(k)}=F_{2n-1} \quad\hbox{for all $1\leq n\leq k+1$}.
\end{equation}
The characteristic polynomial of $P^{(k)}$,
$$
\Psi_k(x) := x^k - 2x^{k-1} - x^{k-2}-\cdots - x - 1,
$$
is irreducible over $\Q$ and has just one root outside the unit circle \cite{BHL2021}. Let
$\alpha := \alpha(k)$ denote that root. It was proved in \cite{BHL2021} that
\begin{equation}\label{eq:alpha}
\varphi^2\left(1-\varphi^{-k}\right)<\alpha<\varphi^2,
\end{equation}
where $\varphi=(1+\sqrt5)/2$ is the golden section. We write $\alpha_1:=\alpha,\alpha_2,\ldots,
\alpha_k$ for the roots of $\Psi_k$; thus $|\alpha_i|<1$ for $2\le i\le k$. A direct computation
gives
\begin{equation}\label{eq:identity}
(x-1)\Psi_k(x)=x^{k+1}-3x^k+x^{k-1}+1,\qquad\hbox{so that}\qquad \alpha^2-3\alpha+1=-\alpha^{1-k}.
\end{equation}
For $x_0=1+\sqrt2$ we have $x_0^2-3x_0+1=1-\sqrt2$, so \eqref{eq:identity} gives
$(x_0-1)\Psi_k(x_0)=1-(\sqrt2-1)x_0^{k-1}=1-x_0^{k-2}\le0$. Since $\Psi_k$ has exactly one
positive root, we get $\alpha\ge1+\sqrt2$, with equality for $k=2$.

For an integer $k\ge2$ we consider the function
\begin{equation}\label{eq:g}
g_{k}(z) = \dfrac{z-1}{(k+1)z^2-3kz+k-1}=\frac{z-1}{k(z^2-3z+1)+z^2-1} \qquad (z \in
\mathbb{C}),
\end{equation}
and we put $g:=g_k(\alpha)$.

\begin{lemma}[{Bravo, Herrera and Luca \cite{BHL2021}}]\label{lem:BHL}
Let $k\geq 2$. Then:
\begin{enumerate}[label=\rm(\alph*)]
\item $\alpha^{n-2}\leq P_n^{(k)}\leq \alpha^{n-1}$ for all $n\geq 1$;
\item $P_n^{(k)} = \sum_{i=1}^{k}g_k(\alpha_i)\alpha_i^n$ for all $n\ge 2-k$;
\item $e_n:=P_n^{(k)} - g\,\alpha^{n}$ satisfies $|e_n|<1/2$ for all $n \geq 2 - k$;
\item $0.276< g<0.5$.
\end{enumerate}
\end{lemma}

\begin{lemma}[{Bravo and Herrera \cite{BH2020}; see also \cite{BHL2021b}}]\label{lem:gconj}
For every $k\ge2$ and $2\leq i\leq k$ we have $$|g_{k}(\alpha_i)|<1.$$
\end{lemma}

We shall also use the following elementary facts. For $n\ge2$ we have $P_n\ge2P_{n-1}$, since
$P_n=2P_{n-1}+(\hbox{non-negative terms})$; hence, for $1\le r\le n$,
\begin{equation}\label{eq:two}
P_r\le2^{-(n-r)}P_n,
\end{equation}
and, for $n\ge2$, $P_{n+1}\ge2P_n+P_{n-1}>2P_n$. The sequence $(P_n)_{n\ge1}$ is strictly
increasing, with $P_1=1$, $P_2=2$ and $P_3=5$ for every $k\ge2$.

\begin{lemma}\label{lem:pow2}
Let $k\ge2$ and $n\ge1$. If $P_n^{(k)}=2^a$ with $a\ge1$, then $n=2$ and $a=1$.
\end{lemma}
\begin{proof}
If $P_n^{(k)}=2^a$, then $|P_n^{(k)}-2^a|=0<2^{a/2}$. Bachabi and Togb\'e \cite[Theorem~1.2]{BT2024}
found all the solutions of $|P_n^{(k)}-2^a|<2^{a/2}$ in positive integers $k\ge2$, $n$, $a$. In
their list the values of $P_n^{(k)}$ are $1,2,5,13,34,29,70,33,4116$, and the only power of $2$
with positive exponent among them is $P_2^{(k)}=2$.
\end{proof}

\begin{lemma}[{Normenyo, Rihane and Togb\'e \cite{NRT2021}}]\label{lem:fermat}
Let $k\ge2$, $m\ge1$ and $a\ge1$. If $P_m^{(k)}=2^a+1$, then $(m,a)=(3,2)$, or $k=3$ and
$(m,a)=(5,5)$.
\end{lemma}

\subsection{Approximation for large $k$}

For large $k$ the dominant root $\alpha$ is very close to $\varphi^2$, and $P_n^{(k)}$ remains very
close to $F_{2n-1}$ also for $n>k+1$, as long as $n$ is small compared with
$\varphi^{2k}$. We put
$$
c:=\frac1{\varphi^2+1}=\frac1{\sqrt5\,\varphi}=\frac{5-\sqrt5}{10}=0.2763\ldots
$$
For instance $F_{2n-1}=c\varphi^{2n}+\varphi^{-2n+1}/\sqrt5$, since $c\varphi=1/\sqrt5$.

\begin{lemma}\label{lem:approx}
Let $k\ge30$ and write $\theta:=\varphi^2-\alpha$, $\delta:=\theta/\varphi^2$ and $g=c(1+\eta)$.
Then
$$
0<\theta<1.1711\,\varphi^{-2k},\qquad 0<\delta<0.4474\,\varphi^{-2k},\qquad
|\eta|\le(0.45k+1.07)\,\varphi^{-2k}.
$$
Consequently, for every $j\ge0$,
\begin{equation}\label{eq:zeta}
g\alpha^j=c\,\varphi^{2j}(1+\zeta_j)\qquad\hbox{with}\qquad
|\zeta_j|\le\big(0.45(k+j)+1.07\big)\varphi^{-2k}.
\end{equation}
\end{lemma}
\begin{proof}
Since $x^2-3x+1=(x-\varphi^2)(x-\varphi^{-2})$, \eqref{eq:identity} gives
$\theta=\alpha^{1-k}/(\alpha-\varphi^{-2})>0$. By \eqref{eq:alpha},
$$\alpha^{1-k}<\varphi^{2-2k}(1-\varphi^{-k})^{1-k}$$ and
$$\alpha-\varphi^{-2}>\varphi^2-\varphi^{-2}-\varphi^{2-k}=\sqrt5-\varphi^{2-k}.$$ For $k\ge30$
this gives $\theta<1.00002\cdot(\varphi^2/\sqrt5)\,\varphi^{-2k}<1.1711\,\varphi^{-2k}$, and
$\delta=\theta/\varphi^2<0.4474\,\varphi^{-2k}$. Next, since $\varphi^4-3\varphi^2+1=0$, we have
$\alpha^2-3\alpha+1=-\theta(\alpha+\varphi^2-3)$, $\alpha^2-1=\varphi^4-1-\theta(\alpha+\varphi^2)$
and $\alpha-1=\varphi^2-1-\theta$. Since $c=(\varphi^2-1)/(\varphi^4-1)$, \eqref{eq:g} gives
$$
\frac gc=\frac{1-u}{1-v},\qquad u=\frac{\theta}{\varphi^2-1},\quad
v=\frac{\theta\,w}{\varphi^4-1},\quad w=(\alpha+\varphi^2)+k(\alpha+\varphi^2-3)<5.3+2.24k.
$$
Here $0<u<v$ and $v<(0.4482k+1.0604)\varphi^{-2k}<10^{-11}$, so that
$$|\eta|=(v-u)/(1-v)\le v/(1-v)\le(0.45k+1.07)\varphi^{-2k}.$$ Finally,
$g\alpha^j=c(1+\eta)\varphi^{2j}(1-\delta)^j$, so $\zeta_j=(1+\eta)(1-\delta)^j-1$. Since
$1-j\delta\le(1-\delta)^j\le1$, we get $-|\eta|-j\delta\le\zeta_j\le|\eta|$, hence
$|\zeta_j|\le|\eta|+0.4474\,j\varphi^{-2k}$, which gives \eqref{eq:zeta}.
\end{proof}

\subsection{Notations and terminology from algebraic number theory}

Let $\eta$ be an algebraic number of degree $d$ with minimal primitive polynomial over the
integers
$$
a_0x^{d}+ a_1x^{d-1}+\cdots+a_d = a_0\prod_{i=1}^{d}(x-\eta^{(i)}),
$$
where the leading coefficient $a_0$ is positive and the $\eta^{(i)}$'s are the conjugates of
$\eta$. Then the \textit{logarithmic height} of $\eta$ is given by
$$
h(\eta):=\dfrac{1}{d}\left( \log a_0 + \sum_{i=1}^{d}\log\left(\max\{|\eta^{(i)}|,
1\}\right)\right).
$$
In particular, if $\eta=p/q$ is a rational number with $\gcd (p,q)=1$ and $q>0$, then
$h(\eta)=\log\max\{|p|, q\}$. The following properties of $h(\cdot)$ will be used without
reference:
\begin{eqnarray*}
h(\eta\pm \gamma) &\leq& h(\eta) +h(\gamma) +\log 2,\\
h(\eta\gamma^{\pm 1})&\leq & h(\eta) + h(\gamma),\\
h(\eta^{s}) &=& |s|h(\eta) \qquad (s\in\mathbb{Z}).
\end{eqnarray*}
Since $\Psi_k$ is monic, irreducible, and has exactly one root of modulus greater than $1$, we
have $h(\alpha)=(\log\alpha)/k$. For $g$ we use the bound of Bravo, Herrera and Luca
\cite[Lemmas~3 and 6]{BHL2021b}:
\begin{equation}\label{eq:hg}
h(g)<H(k):=4\log\varphi+\log(k+1)\qquad(k\ge2).
\end{equation}
The next lemma is needed for the case $r=n$ of \eqref{eq:main}.

\begin{lemma}\label{lem:disc}
Let $k\ge2$. Then $g=\alpha^{k-2}/\Psi_k'(\alpha)$, so that $|N_{\Q(\alpha)/\Q}(g)|=1/|\Delta_k|$,
where $\Delta_k$ is the discriminant of $\Psi_k$. Moreover $|\Delta_k|>2^k$.
\end{lemma}
\begin{proof}
Put $Q(z):=(z-1)\Psi_k(z)=z^{k+1}-3z^k+z^{k-1}+1$ and $q(z):=(k+1)z^2-3kz+k-1$. Then
$Q'(z)=z^{k-2}q(z)$, and differentiating $Q=(z-1)\Psi_k$ at $z=\alpha_i$ gives
$(\alpha_i-1)\Psi_k'(\alpha_i)=\alpha_i^{k-2}q(\alpha_i)$. By \eqref{eq:g}, $g=(\alpha-1)/q(\alpha)
=\alpha^{k-2}/\Psi_k'(\alpha)$. Since $|N(\alpha)|=1$ and
$|N(\Psi_k'(\alpha))|=|\Delta_k|$, we get $|N(g)|=1/|\Delta_k|$. Next, let $\rho_1>\rho_2$ be the
roots of $q$. Taking the product of $(\alpha_i-1)\Psi_k'(\alpha_i)=\alpha_i^{k-2}q(\alpha_i)$ over
all $i$ and using $\prod_i|\alpha_i|=1$, $\prod_i|\alpha_i-1|=|\Psi_k(1)|=k$ and
$$\prod_iq(\alpha_i)=(k+1)^k\Psi_k(\rho_1)\Psi_k(\rho_2),$$ we obtain
$$
|\Delta_k|=\frac{(k+1)^k}{k}\,|\Psi_k(\rho_1)|\,|\Psi_k(\rho_2)|.
$$
Since $q(\varphi^{\pm2})=\varphi^{\pm4}-1$, we have $0<\rho_2<\varphi^{-2}<\rho_1<\varphi^2$. As
$\rho_2^2-3\rho_2+1>0$, we get $Q(\rho_2)>1$ and $|\Psi_k(\rho_2)|=|Q(\rho_2)|/(1-\rho_2)>1$. Let
$k\ge5$. Then $$\rho_1=\big(3k+\sqrt{5k^2+4}\big)/(2k+2)>\varphi^2k/(k+1)>2.18$$ and
$\varphi^2-\rho_1>2/(k+1)$, so that
$$\rho_1^{k-1}(\varphi^2-\rho_1)(\rho_1-\varphi^{-2})>2^{k-1}\cdot\frac2{k+1}\cdot1.7>9;$$ since
$\rho_1^2-3\rho_1+1=-(\varphi^2-\rho_1)(\rho_1-\varphi^{-2})$, this gives $|Q(\rho_1)|>8$ and
$|\Psi_k(\rho_1)|=|Q(\rho_1)|/(\rho_1-1)>4$. Hence $|\Delta_k|>4(k+1)^k/k>2^k$. For $k=2,3,4$ we
have $|\Delta_k|=8,87,1423$.
\end{proof}

In $\Q(\sqrt5)$ we shall use $h(2)=\log2$, $h(\sqrt5)=\frac12\log5$, $h(\varphi)=\frac12\log\varphi$,
$h(2/\sqrt5)=\frac12\log5$, and
$h(2c)=\frac12\log\big(5(1+1/\sqrt5)\big)<0.98953$, since $2c=1-1/\sqrt5$ has minimal
polynomial $5y^2-10y+4$.

\subsection{Linear forms in logarithms and the reduction lemma}

We use the following result of Matveev \cite{Matveev2000}.

\begin{theorem}[Matveev]\label{thm:matveev}
Let $\gamma_1,\ldots,\gamma_t$ be positive real algebraic numbers in a real algebraic number
field $\mathbb{K}$ of degree $D$, $b_1,\ldots,b_t$ be nonzero integers, and assume that
$$
\Lambda:=\gamma_1^{b_1}\cdots\gamma_t^{b_t} - 1
$$
is nonzero. Then
$$
\log |\Lambda| > -1.4\times 30^{t+3}\times t^{4.5}\times D^{2}(1+\log D)(1+\log B)A_1\cdots A_t,
$$
where $B\geq\max\{|b_1|, \ldots, |b_t|\}$ and $A_i \geq \max\{Dh(\gamma_i), |\log\gamma_i|,
0.16\}$ for all $i=1,\ldots,t$.
\end{theorem}

We write $C(D,t):=1.4\times30^{t+3}\times t^{4.5}\times D^2(1+\log D)$. For a real number $r$ we denote by $\| r\|$ the distance from $r$ to the nearest integer. The
following is the reduction lemma of Dujella and Peth\H{o} \cite[Lemma~5(a)]{DP1998}; we state it
with a real exponent $w$, the proof being the same.

\begin{lemma}[Dujella and Peth\H{o}]\label{lem:dp}
Let $M$ be a positive integer, let $\tau$ be an irrational number and let $p/q$ be a convergent
of the continued fraction of $\tau$ with $q>6M$. Let $A,B,\mu$ be real numbers with $A>0$ and
$B>1$, and put $\varepsilon:=\|\mu q\|-M\|\tau q\|$. If $\varepsilon>0$, then there is no
solution of
$$
0<|u\tau-v+\mu|<AB^{-w}
$$
in integers $u,v$ and a real number $w$ with $1\le u\le M$ and $w\ge\log(Aq/\varepsilon)/\log B$.
\end{lemma}

In practice we take the first convergent with $q>6M$, and the next ones if $\varepsilon\le0$.
Lemma~\ref{lem:dp} cannot be used when $\mu$ is (almost) a linear combination of $1$ and $\tau$;
in that case we use the following property of continued fractions
\cite[Theorem~8.2.4]{ME2005}.

\begin{lemma}\label{lem:leg}
Let $\tau$ be irrational, with continued fraction $[a_0;a_1,a_2,\ldots]$ and convergents
$p_i/q_i$. Let $M$ be a positive integer and $N$ such that $q_N\le M<q_{N+1}$, and put
$a(M):=\max\{a_i:0\le i\le N+1\}$. Then $|u\tau-v|>1/((a(M)+2)u)$ for all integers $u,v$ with
$0<u<M$.
\end{lemma}


\section{Elementary cases and the size of $m$}\label{sec:elem}

Let $(k,n,r,m,x)$ be a solution of \eqref{eq:main} in positive integers with $k\ge2$ and $n\ge r$.
We write $d:=n-r\ge0$.

\begin{lemma}[The case $x=1$]\label{lem:x1}
If $x=1$, then $(n,r,m)=(1,1,2)$.
\end{lemma}
\begin{proof}
If $n=1$, then $r=1$ and $P_m=2$, so $m=2$. If $n\ge2$, then
$P_n<P_n+P_r\le2P_n<P_{n+1}$, and $P_m$ cannot lie strictly between two consecutive terms.
\end{proof}

\begin{lemma}[The case $n\le2$]\label{lem:n2}
If $n\le2$ and $x\ge2$, then $(n,r,m)=(1,1,2)$, or $(n,r,m,x)=(2,1,3,2)$, or $k=3$ and
$(n,r,m,x)=(2,1,5,5)$.
\end{lemma}
\begin{proof}
If $n=1$, then $r=1$ and $m=2$. Let $n=2$. If $r=2$, then $P_m=2^{x+1}$, which is impossible by
Lemma~\ref{lem:pow2} since $x\ge1$. If $r=1$, then $P_m=2^x+1$, and Lemma~\ref{lem:fermat} gives
$(m,x)=(3,2)$, or $k=3$ and $(m,x)=(5,5)$.
\end{proof}

\begin{lemma}[The case $k=2$, $x=2$, $d=1$]\label{lem:pell}
If $k=2$, $x=2$ and $r=n-1\ge1$, then $m=2n-1$.
\end{lemma}
\begin{proof}
By \eqref{eq:pellid}, $P_m=P_{n-1}^2+P_n^2=P_{2n-1}$.
\end{proof}

The values listed in Theorem~\ref{thm:main} are solutions; for (c) this is \eqref{eq:pellid}. By
the three lemmas, it remains to show that \eqref{eq:main} has no solution with
\begin{equation}\label{eq:remaining}
n\ge3,\qquad x\ge2,\qquad 1\le r\le n,\qquad (k,x,d)\ne(2,2,1).
\end{equation}

\begin{lemma}[Size of $m$]\label{lem:size}
If $n\ge3$ and $x\ge2$, then $(n-2)x+2\le m\le(n-1)x+2$. In particular $m\ge x+2$, and the
integer
$$
j:=nx-m
$$
satisfies $x-2\le j\le2x-2$.
\end{lemma}
\begin{proof}
By Lemma~\ref{lem:BHL}(a), $$\alpha^{(n-2)x}\le P_n^x<P_m\le\alpha^{m-1},$$ so $m-1>(n-2)x$. Also
$\alpha^{m-2}\le P_m\le2P_n^x\le2\alpha^{(n-1)x}$, so $m-2\le(n-1)x+\log2/\log\alpha<(n-1)x+1$.
\end{proof}
\section{Proof of Theorem~\ref{thm:fib}}\label{sec:fib}

We begin with the case in which all the indices are at most $k+1$, which, as explained in the
introduction, reduces to \eqref{eq:fib}. The argument has the same structure as the general one,
but everything takes place in $\Q(\sqrt5)$, the heights are small and the constants are as good as
they can be. This is also where we first met the difficulties mentioned in the introduction, so
we give more details here than in the following sections.

In this section $(n,r,m,x)$ is a solution of \eqref{eq:fib} with $n\ge3$, $1\le r\le n$ and
$x\ge2$, and $d:=n-r$. For $s\ge1$ we have
$$F_{2s-1}=\frac{\varphi^{2s-1}}{\sqrt5}(1+a_s)$$ with $a_s:=\varphi^{-4s+2}$, and
$\varphi^{t-2}\le F_t\le\varphi^{t-1}$ for $t\ge1$. Since $F_{t+2}\ge2F_t$, we have
$$F_{2r-1}\le2^{-d}F_{2n-1}.$$

We first give a bound on $m$. From $\varphi^{(2n-3)x}\le F_{2n-1}^x<F_{2m-1}\le\varphi^{2m-2}$ and
$$\varphi^{2m-3}\le F_{2m-1}\le2F_{2n-1}^x\le2\varphi^{(2n-2)x}$$ we get, 
$$
e:=2m-1-(2n-1)x,\qquad 2-2x\le e\le3-x,\qquad\hbox{and}\qquad m\le(n-1)x+2 .
$$
Note also that $e\equiv x+1\pmod2$.

 Each of the following forms has a precise role. $\Phi_A$ compares
$F_{2m-1}$ with $F_{2n-1}^x$,  $\Phi_B$ compares $F_{2m-1}$ with the main part $(\varphi^{2n-1}/\sqrt5)^x$ of
$F_{2n-1}^x$ and $\Phi_C$ bounds $n$. The forms with
index $0$ play the same roles when $r=n$. We use
\begin{align*}
\Phi_A&:=\frac{\varphi^{2m-1}}{\sqrt5}\,F_{2n-1}^{-x}-1, &
\Phi_{A0}&:=\frac{\varphi^{2m-1}}{2\sqrt5}\,F_{2n-1}^{-x}-1,\\
\Phi_B&:=5^{(x-1)/2}\varphi^{e}-1, &
\Phi_{B0}&:=\tfrac12\,5^{(x-1)/2}\varphi^{e}-1,\\
\Phi_C&:=5^{(x-1)/2}\varphi^{e}\,(1+\varphi^{-2dx})^{-1}-1, &&
\end{align*}
the forms with index $0$ being used when $d=0$, the others when $d\ge1$. Consider the conditions
\begin{align*}
20x&\le\varphi^{4n-2},\tag{$\ast$}\label{eq:starF}\\
20x&\le\varphi^{4r-2}.\tag{$\ast\ast$}\label{eq:starFr}
\end{align*}

\begin{lemma}\label{lem:Phi}
\begin{enumerate}[label=\rm(\alph*)]
\item If $d\ge1$, then $|\Phi_A|<1.1\cdot2^{-dx}$; if $d=0$, then $|\Phi_{A0}|<0.1\cdot25^{-x}$.
\item Under \eqref{eq:starF}: if $d\ge1$, then $|\Phi_B|<2.1\max\{2^{-dx},xa_n\}$; if $d=0$, then
$|\Phi_{B0}|<1.03\,xa_n$.
\item Under \eqref{eq:starFr}, if $d\ge1$, then $|\Phi_C|<2.1\,xa_n$.
\item The five forms are non-zero.
\end{enumerate}
\end{lemma}
\begin{proof}
(a) Since $\varphi^{2m-1}/\sqrt5=F_{2m-1}-\varphi^{1-2m}/\sqrt5$, we have
$$
\Phi_A=\Big(\frac{F_{2r-1}}{F_{2n-1}}\Big)^x-\frac{\varphi^{1-2m}}{\sqrt5\,F_{2n-1}^x}.
$$
The first term is at most $2^{-dx}$. Since $\varphi^{2m-1}>F_{2m-1}>F_{2n-1}^x$ and
$F_{2n-1}^2\ge2^{n+1}$ for $n\ge3$, the second one is less than
$2^{-dx}\,4^{-x}/\sqrt5\le0.03\cdot2^{-dx}$. If $d=0$, then $F_{2m-1}=2F_{2n-1}^x$, so
$\Phi_{A0}=-\varphi^{1-2m}/(2\sqrt5F_{2n-1}^x)$ and $\varphi^{2m-1}>4.4F_{2n-1}^x$; since
$F_{2n-1}\ge5$, $|\Phi_{A0}|<25^{-x}/19$.

(b) Dividing \eqref{eq:fib} by $(\varphi^{2n-1}/\sqrt5)^x$ we get
$$(\Phi_B+1)(1+a_m)=(1+a_n)^x+\big(F_{2r-1}\sqrt5/\varphi^{2n-1}\big)^x.$$ Under
\eqref{eq:starF}, $xa_n\le0.05$, so $(1+a_n)^x-1\le1.026\,xa_n$ and
$$(F_{2r-1}\sqrt5/\varphi^{2n-1})^x=(F_{2r-1}/F_{2n-1})^x(1+a_n)^x\le1.052\cdot2^{-dx}$$ moreover
$a_m\le\varphi^{-8n+8}$. This gives the bound for $\Phi_B$, and the bound for
$\Phi_{B0}$ is obtained in the same way from $(\Phi_{B0}+1)(1+a_m)=(1+a_n)^x$.

(c) In the same way, $$(\Phi_C+1)(1+\varphi^{-2dx})(1+a_m)-(1+\varphi^{-2dx})=
[(1+a_n)^x-1]+\varphi^{-2dx}[(1+a_r)^x-1].$$Under \eqref{eq:starFr}, the second bracket is at
most $1.026\,xa_r$, and $\varphi^{-2dx}a_r=a_n\varphi^{2d(2-x)}\le a_n$. Hence $|\Phi_C|<(1.026+
1.026+0.01)\,xa_n$.

(d) We use the norm $N$ from $\Q(\sqrt5)$ to $\Q$, with $|N(\varphi)|=1$ and $|N(\sqrt5)|=5$. If
$\Phi_A=0$ or $\Phi_{A0}=0$, then $F_{2n-1}^{2x}=|N(F_{2n-1}^x)|$ would be equal to $1/5$ or
$1/20$. If $\Phi_B=0$ or $\Phi_{B0}=0$, then $5^{x-1}=1$ or $5^{x-1}=4$. If $\Phi_C=0$, put
$t:=dx\ge2$; then $$5^{x-1}=|N(1+\varphi^{-2t})|=(1+\varphi^{-2t})(1+\varphi^{2t})=L_{2t}+2.$$ Now
$L_{2t}+2=L_t^2$ if $t$ is even and $L_{2t}+2=5F_t^2$ if $t$ is odd. Since $5\nmid L_t$, the first
case is impossible for $x\ge2$. In the second, $F_t=5^{(x-2)/2}$ forces $F_t\in\{1,5\}$, that is
$t=1$ or $(t,x)=(5,4)$, which is impossible since $x\mid t$.
\end{proof}

We apply Theorem~\ref{thm:matveev} in $\mathbb K=\Q(\sqrt5)$, $D=2$, with
$h(\sqrt5)=\frac12\log5$, $h(\varphi)=\frac12\log\varphi$, $h(2\sqrt5)=\frac12\log20$ and
$$h(1+\varphi^{-2t})=\frac12\log(L_{2t}+2)\le t\log\varphi+\frac12$$  With $B=2x$ for $\Phi_B,\Phi_{B0},\Phi_C$ and $B=2m$ for
$\Phi_A,\Phi_{A0}$, and $2\log F_{2n-1}\le4(n-1)\log\varphi$, we obtain
{\small\begin{align*}
\log|\Phi_B|&>-K_B(1+\log2x), & K_B&:=C(2,2)\log5\log\varphi<4.038\times10^{9},\\
\log|\Phi_C|&>-K_C(2dx\log\varphi+1)(1+\log2x), & K_C&:=C(2,3)\log5\log\varphi<7.511\times10^{11},\\
\log|\Phi_{B0}|&>-K_{B0}(1+\log2x), & K_{B0}&:=C(2,3)\cdot2\log2\cdot\log5\log\varphi<1.042\times10^{12},\\
\log|\Phi_A|&>-4(n-1)\log\varphi\,K_C(1+\log2m), &&\hbox{(same $K_C$),}\\
\log|\Phi_{A0}|&>-4(n-1)\log\varphi\,K_{A0}(1+\log2m), & K_{A0}&:=C(2,3)\log20\log\varphi<1.399\times10^{12}.
\end{align*}}

Now let $d\ge1$. If \eqref{eq:starF} holds and the minimum in
Lemma~\ref{lem:Phi}(b) is $2^{-dx}$, then $$dx\log2<\log2.1+K_B(1+\log2dx)$$ hence $x\le
dx<1.602\times10^{11}$. Otherwise, either \eqref{eq:starF} fails, or
$$(4n-2)\log\varphi-\log x<\log2.1+K_B(1+\log2x)$$ in both cases $n$ is bounded by an explicit
function of $x$. Inserting this and $m\le(n-1)x+2$ in the bound for
$\Phi_A$, we get $x<2.251\times10^{25}$. Let $d=0$. Then either \eqref{eq:starF} fails or
$$(4n-2)\log\varphi-\log(1.03x)<K_{B0}(1+\log2x),$$ and the bound for $\Phi_{A0}$ gives
$x<2.820\times10^{27}$. Hence in all cases $x<2.82\times10^{27}$.

Now we can make the reduction procedure. One pass, with $x\le M$, consists of the following steps.
\begin{enumerate}[label=\rm(F\arabic*)]
\item ($d\ge1$.) With $\tau_0:=\log\sqrt5/\log\varphi$, Lemma~\ref{lem:Phi}(b) gives
$|(x-1)\tau_0+e|<(4.2/\log\varphi)e^{-w}$ with $w=\min\{dx\log2,(4n-2)\log\varphi-\log x\}$.
Lemma~\ref{lem:leg} with $u=x-1$ gives $w<W$. Hence either $dx<W/\log2$, or $n$ is bounded, or
\eqref{eq:starF} fails and $4n<2+\log(20M)/\log\varphi$.
\item ($d\ge1$, $dx<W/\log2$.) For each such $d$, $x$ and each $e\in[2-2x,3-x]$ with $e\equiv x+1\pmod2$, the number
$\Phi_C$ is explicit and non-zero, and Lemma~\ref{lem:Phi}(c) bounds $n$; if
\eqref{eq:starFr} fails, then $r$ is bounded and so is $n=r+d$.
\item ($d=0$.) Lemma~\ref{lem:dp} applied to $(x-1)\tau_0+e-\log2/\log\varphi$, with $u=x-1$,
$A=2.06M/\log\varphi$, $B=\varphi$ and $w=4n-2$, bounds $n$.
\item For each $3\le n\le N$ (the bound obtained in (F1)--(F3)), Lemma~\ref{lem:dp} applied to
$\Phi_A$ (with $u=2m-1$, $B=2$, $w=dx$) and to $\Phi_{A0}$ (with $B=25$, $w=x$) gives a new bound
for $x$.
\end{enumerate}
In (F4), for large $n$, Lemma~\ref{lem:dp} fails for $\Phi_A$ because 
$$\log F_{2n-1}=(2n-1)\log\varphi-\log\sqrt5+\log(1+a_n),$$ the number
$\mu=-\log\sqrt5/\log F_{2n-1}$ is extremely close to $1-(2n-1)\tau$, where
$\tau=\log\varphi/\log F_{2n-1}$. For these $n$ we write, with $\Gamma_A:=\log(1+\Phi_A)$,
$$
\Gamma_A/\log\varphi=-\big((x-1)\tau_n-(2m-2n)+\delta_n\big),\qquad
\tau_n=\frac{\log F_{2n-1}}{\log\varphi},\quad \delta_n=\frac{\log(1+a_n)}{\log\varphi},
$$
and use Lemma~\ref{lem:leg} with $u=x-1$, in all
these cases $\delta_n\le1/(2(a(M)+2)M)$. Table~\ref{tab:fibred} shows the passes.

\begin{table}[ht]
\centering
\caption{Reduction for \eqref{eq:fib}. Here $M$ is the bound for $x$ used in the pass, $W$ comes
from (F1), $N$ is the bound for $n$ and the last two columns give the bounds for $x$ obtained
from $\Phi_A$ and $\Phi_{A0}$.}\label{tab:fibred}
\begin{tabular}{@{}rrrrrrr@{}}
\toprule
pass & $M$ & $W$ & $dx<W/\log2$ & $N$ & $\Phi_A$ & $\Phi_{A0}$\\
\midrule
1 & $2.82\times10^{27}$ & $68.81$ & $99$ & $69$ & $160$ & $22$\\
2 & $161$ & $10.30$ & $14$ & $9$ & $24$ & $4$\\
3 & $25$ & $8.47$ & $12$ & $8$ & $23$ & $3$\\
4 & $24$ & $8.43$ & $12$ & $8$ & $22$ & $3$\\
5 & $23$ & $8.39$ & $12$ & $8$ & $22$ & $3$\\
\bottomrule
\end{tabular}
\end{table}

As usual, the first pass does almost all the work. In fact,  the bound for $x$ falls from $10^{27}$ to $160$,
Thus $n\le8$ and $x\le22$. A direct search
with $3\le n\le40$, $1\le r\le n$ and $2\le x\le40$ shows that \eqref{eq:fib} has no solution.
This proves Theorem~\ref{thm:fib}.\qed
\section{Linear forms in $\Q(\alpha)$ and a bound for $x$}\label{sec:generic}

We now leave the Fibonacci numbers and work with $P^{(k)}$ for an arbitrary $k$. The forms are those
of Section~\ref{sec:fib}, with $g$ and $\alpha$ in the place of $1/\sqrt5$ and $\varphi^2$. We now work in a field of degree $k$, and that the height of $g$ grows like
$\log k$, so that the bound we obtain for $x$ depends on $k$.

In this section $k\ge2$ and $(k,n,r,m,x)$ is a solution of \eqref{eq:main} satisfying
\eqref{eq:remaining}, @$d=n-r$ and $j=nx-m$ are as in Section~\ref{sec:elem}. We consider
\begin{align*}
\Lambda_A&:=g\,\alpha^{m}P_n^{-x}-1, & \Lambda_{A0}&:=\tfrac12\,g\,\alpha^{m}P_n^{-x}-1,\\
\Lambda_B&:=g^{1-x}\alpha^{-j}-1, & \Lambda_{B0}&:=\tfrac12\,g^{1-x}\alpha^{-j}-1,\\
\Lambda_C&:=g^{1-x}\alpha^{-j}(1+\alpha^{-dx})^{-1}-1, &&
\end{align*}
(the forms with index $0$ for $d=0$, the others for $d\ge1$) and the conditions
\begin{align*}
36.4\,x&\le\alpha^{n-1},\tag{$a\ast$}\label{eq:star}\\
36.4\,x&\le\alpha^{r-1}.\tag{$b\ast$}\label{eq:starr}
\end{align*}

\begin{lemma}\label{lem:Lambda}
\begin{enumerate}[label=\rm(\alph*)]
\item If $d\ge1$, then $0<|\Lambda_A|<1.5\cdot2^{-dx}$; if $d=0$, then
$0<|\Lambda_{A0}|<0.25\cdot5^{-x}$.
\item Under \eqref{eq:star}: if $d\ge1$, then $0<|\Lambda_B|<3\max\{2^{-dx},x\alpha^{-n}\}$; if
$d=0$, then $0<|\Lambda_{B0}|<1.9\,x\alpha^{-n}$.
\item Under \eqref{eq:starr}, if $d\ge1$, then $|\Lambda_C|<2.7\,x\alpha^{-n}$.
\end{enumerate}
\end{lemma}
\begin{proof}
(a) Since $g\alpha^m=P_m-e_m=P_n^x+P_r^x-e_m$, we have
$\Lambda_A=(P_r/P_n)^x-e_mP_n^{-x}$, and \eqref{eq:two} gives $|\Lambda_A|<2^{-dx}+\frac12
2^{-dx}$. If $d=0$, then $$\Lambda_{A0}=-e_m/(2P_n^x)$$ and $P_n\ge5$. If $\Lambda_A=0$, then
$g\alpha^m=P_n^x$. Then, applying an automorphism of the Galois closure of $\Q(\alpha)$ which sends
$\alpha$ to $\alpha_2$ and using Lemma~\ref{lem:gconj} we get $P_n^x<1$, a contradiction.  The same
argument works for $\Lambda_{A0}$.

(b) Put $U:=g\alpha^n$ and $P_n=U(1+\varepsilon)$. By Lemma~\ref{lem:BHL},
$|\varepsilon|<1.812\,\alpha^{-n}$, and \eqref{eq:star} gives $x|\varepsilon|\le0.05$ and
$n\ge6$. Dividing \eqref{eq:main} by $U^x$,
$$
\Lambda_B=\big[(1+\varepsilon)^x-1\big]+(P_r/P_n)^x(1+\varepsilon)^x-e_mU^{-x}.
$$
Here $|(1+\varepsilon)^x-1|\le1.026\,x|\varepsilon|<1.86\,x\alpha^{-n}$, the second term is less
than $1.052\cdot2^{-dx}$, and $|e_m|U^{-x}<0.04\,x\alpha^{-n}$. This gives the bound for
$\Lambda_B$. Hence, for $d=0$, $$\Lambda_{B0}=[(1+\varepsilon)^x-1]-e_m/(2U^x).$$ If $\Lambda_B=0$, then
$g^{x-1}\alpha^j=1$, and conjugating as in (a) we get $|g_k(\alpha_2)|^{x-1}|\alpha_2|^j=1$, which
is impossible since $x\ge2$ and $j\ge0$. If $\Lambda_{B0}=0$, then $2g^{x-1}\alpha^j=1$, and
taking norms, $2^k|N(g)|^{x-1}=1$, that is, $|\Delta_k|^{x-1}=2^k$. This contradicts
Lemma~\ref{lem:disc}.

(c) Put $V:=g\alpha^r$ and $P_r=V(1+\varepsilon_r)$, so that $V/U=\alpha^{-d}$. Then
$$
\Lambda_C\,(1+\alpha^{-dx})=\big[(1+\varepsilon)^x-1\big]+\alpha^{-dx}\big[(1+\varepsilon_r)^x-1\big]-e_mU^{-x}.
$$
Under \eqref{eq:starr}, $\alpha^{-dx}|(1+\varepsilon_r)^x-1|<1.86\,x\alpha^{-r-dx}=1.86\,x\alpha^{-n}
\alpha^{-d(x-1)}<0.78\,x\alpha^{-n}$, and the claim follows.
\end{proof}

\begin{lemma}\label{lem:zero}
\begin{enumerate}[label=\rm(\alph*)]
\item For $k=2$ and $(d,x,j)=(1,2,1)$ we have $\Lambda_C=0$.
\item In all the cases met in step (S2) of Section~\ref{sec:smallk}, that is, for $2\le k\le K_0$
and the triples $(d,x,j)$ used there, $\Lambda_C\ne0$ except in the case (a), and
$|\Lambda_C|>8.9\times10^{-6}$.
\end{enumerate}
\end{lemma}
\begin{proof}
(a) For $k=2$ we have $\alpha=1+\sqrt2$ and $g=1/(2\sqrt2)$, so that
$$g^{-1}\alpha^{-1}=2\sqrt2\,\alpha^{-1}=1+\alpha^{-2}$$ because $\alpha+\alpha^{-1}=2\sqrt2$. This
is \eqref{eq:pellid} seen in $\Q(\sqrt2)$. (b) This is a direct computation
(Appendix~\ref{app:comp}).
\end{proof}

The zero in (a) is excluded by \eqref{eq:remaining}. We use $\Lambda_C$ only for explicit triples
$(d,x,j)$, so that a general non-vanishing result is not needed.

We can now apply Matveev's theorem. Let $\mathbb K=\Q(\alpha)$, $D=k$. For $\Lambda_B$ we take $t=2$,
$(\gamma_1,\gamma_2)=(g,\alpha)$, $(b_1,b_2)=(1-x,-j)$, $B=2x$, $A_1=kH(k)$, $A_2=\log\alpha$. For
$\Lambda_{B0}$ we add $\gamma_3=2$, with $A_3=k\log2$. For $\Lambda_A$ we take $t=3$,
$(\gamma_1,\gamma_2,\gamma_3)=(g,\alpha,P_n)$, $(b_1,b_2,b_3)=(1,m,-x)$, $B=m$ (since $m>x$),
$A_1=kH(k)$, $A_2=\log\alpha$ and $A_3=k\log P_n\le k(n-1)\log\alpha$; for $\Lambda_{A0}$,
$\gamma_1=g/2$ and $A_1=k(H(k)+\log2)$. Using $\log\alpha<2\log\varphi$ we get, with
$\ell:=2\log\varphi$,
\begin{align}
\min\{dx\log2,\ n\log\alpha-\log x\}-\log3&<E_B(k)(1+\log2x),\label{eq:MB}\\
n\log\alpha-\log(1.9x)&<E_{B0}(k)(1+\log2x),\label{eq:MB0}\\
dx\log2-\log1.5&<E_A(k)\,(n-1)(1+\log m),\label{eq:MA}\\
x\log5+\log4&<E_{A0}(k)\,(n-1)(1+\log m),\label{eq:MA0}
\end{align}
the first two under \eqref{eq:star}, where
\begin{align*}
E_B(k)&:=C(k,2)\,kH(k)\,\ell, & E_{B0}(k)&:=C(k,3)\,k\log2\cdot kH(k)\,\ell,\\
E_A(k)&:=C(k,3)\,k^2H(k)\,\ell^2, & E_{A0}(k)&:=C(k,3)\,k^2(H(k)+\log2)\,\ell^2.
\end{align*}

\begin{proposition}\label{prop:X}
Let $X(k)$ be the largest of the following three numbers: 
\begin{itemize}
\item the supremum $Y(k)$ of the real $t\ge2$
with $t\log2<\log3+E_B(k)(1+\log2t)$; 
\item the supremum of the real $x\ge2$ with
$$x\log2<\log1.5+E_A(k)(N_s(x)-1)\big(1+\log((N_s(x)-1)x+2)\big)$$ where
$N_s(x):=1+\big(\log(36.4x)+\log3+E_B(k)(1+\log2x)\big)/\log(1+\sqrt2)$ 
\item the supremum of the
real $x\ge2$ with $$x\log5+\log4<E_{A0}(k)(N_0(x)-1)\big(1+\log((N_0(x)-1)x+2)\big)$$
 where
$N_0(x):=1+\big(\log(36.4x)+\log1.9+E_{B0}(k)(1+\log2x)\big)/\log(1+\sqrt2)$. 
\end{itemize}
Then every
solution of \eqref{eq:main} satisfying \eqref{eq:remaining} has $x<X(k)$. For instance
$X(2)<5.35\times10^{29}$, $X(30)<4.97\times10^{40}$ and $X(604)<1.44\times10^{52}$.
\end{proposition}
\begin{proof}
Let $d\ge1$. If \eqref{eq:star} holds and the minimum in \eqref{eq:MB} is $dx\log2$, then
$t:=dx\ge x$ satisfies the first inequality. Otherwise $n<N_s(x)$, and \eqref{eq:MA} with
$dx\ge x$ and $m\le(n-1)x+2$ gives the second inequality. If $d=0$, then either \eqref{eq:star}
fails or \eqref{eq:MB0} holds, so $n<N_0(x)$, and \eqref{eq:MA0} gives the third one.
\end{proof}
\section{Large values of $k$}\label{sec:largek}

In this section, we
replace the linear forms of Section~\ref{sec:generic} by those of Section~\ref{sec:fib}, whose
coefficients do not depend on $k$. The cost is an error term of about $k\varphi^{-2k}$. This error
does no harm once $k$ is larger than a few hundred.

Here, we assume that $k\ge30$. By Lemma~\ref{lem:approx}, we can replace $g$ by
$c=1/(\sqrt5\varphi)$ and $\alpha$ by $\varphi^2$. We have
$c^{1-x}\varphi^{-2j}=5^{(x-1)/2}\varphi^{x-1-2j}$ and $x-1-2j=2m-1-(2n-1)x=e$. So the forms
$\Lambda_B$, $\Lambda_{B0}$ and $\Lambda_C$ become the forms $\Phi_B$, $\Phi_{B0}$ and $\Phi_C$ of
Section~\ref{sec:fib}. Here $e$ now satisfies $3-3x\le e\le3-x$, by Lemma~\ref{lem:size}.
When $n\le k+1$, we have $P_n=F_{2n-1}$ by \eqref{eq:kilic}. Since
$c\varphi^{2m}=\varphi^{2m-1}/\sqrt5$, the forms $\Lambda_A$ and $\Lambda_{A0}$ then become
$\Phi_A$ and $\Phi_{A0}$. By Lemma~\ref{lem:Phi}(d), all these forms are non-zero. Matveev's
theorem therefore gives the bounds of Section~\ref{sec:fib}, with $1+\log3x$ instead of
$1+\log2x$. We set
$$
E_k(x):=1.011\,x\,(0.45k+2.02)\,\varphi^{-2k},\qquad\lambda_*:=-\log\big(2.022E_k(x)\big),
$$
$$
\lambda_\zeta:=-\log\big(2.1(0.45(k+m)+1.07)\varphi^{-2k}\big).
$$

\begin{lemma}\label{lem:transfer}
Let $k\ge30$, and let $(k,n,r,m,x)$ be a solution of \eqref{eq:main} which satisfies
\eqref{eq:remaining} and $2.022\,E_k(x)\le0.01$.
\begin{enumerate}[label=\rm(\alph*)]
\item Under \eqref{eq:star}: if $d\ge1$, then
$|\Phi_B|<e^{-\lambda_B}$ with $$\lambda_B:=\min\{dx\log2-\log6.07,\ n\log\alpha-\log(6.07x),\
\lambda_*\}$$ if $d=0$, then $|\Phi_{B0}|<e^{-\lambda_{B0}}$ with
$\lambda_{B0}:=\min\{n\log\alpha-\log(3.84x),\ \lambda_*\}$.
\item Under \eqref{eq:starr}, if $d\ge1$, then $|\Phi_C|<e^{-\lambda_C}$ with
$\lambda_C:=\min\{n\log\alpha-\log(5.46x),\ \lambda_*\}$.
\item If $n\le k+1$ and $|\zeta_m|\le0.01$, then $|\Phi_A|<e^{-\lambda_A}$ with
$\lambda_A:=\min\{dx\log2-\log3.04,\ \lambda_\zeta\}$ if $d\ge1$, and $|\Phi_{A0}|<e^{-\lambda_{A0}}$
with $\lambda_{A0}:=\min\{x\log5-\log0.506,\ \lambda_\zeta\}$ if $d=0$.
\end{enumerate}
\end{lemma}
\begin{proof}
By Lemma~\ref{lem:approx}, $$g^{1-x}\alpha^{-j}=(\Phi_B+1)\rho$$ with
$\rho=(1+\eta)^{1-x}(1-\delta)^{-j}$, and, since $j\le2x-2$,
$$|\log\rho|\le1.0001\big((x-1)|\eta|+j\delta\big)\le1.0001\,x(0.45k+1.965)\varphi^{-2k}$$ so that
$|\rho-1|\le E_k(x)$. Hence $$|\Phi_B|=|\Lambda_B+1-\rho|/\rho\le2.0204\max\{|\Lambda_B|,E_k(x)\}$$
and Lemma~\ref{lem:Lambda} gives (a). For
$\Phi_C$ there is the additional factor $(1+\varphi^{-2dx})/(1+\alpha^{-dx})$, which differs from
$1$ by less than $1.02\,dx\delta\varphi^{-2dx}<0.2\,\varphi^{-2k}$ and is absorbed in $E_k(x)$. For (c), $$\Lambda_A+1=(\Phi_A+1)(1+\zeta_m)$$ by \eqref{eq:zeta}, so
$|\Phi_A|\le(|\Lambda_A|+|\zeta_m|)/(1-|\zeta_m|)\le2.0202\max\{1.5\cdot2^{-dx},|\zeta_m|\}$, and
similarly for $\Phi_{A0}$.
\end{proof}

\begin{lemma}\label{lem:Esmall}
For every $k\ge120$ we have $2.022\,E_k(X(k))<0.01$.
\end{lemma}
\begin{proof}
Let $k\ge120$. Since $(1+\log k)(\log k+2.627)\le k^{0.8}$ for $k\ge120$ (at $k=120$ the two sides
are $<42.92$ and $>46.06$, and their ratio increases with $k$), the definitions give
\[
\max\{E_B(k),E_{B0}(k),E_A(k),E_{A0}(k)\}\le M:=1.33\cdot10^{11}k^{4.8},\qquad M\ge1.27\cdot10^{21}.
\]
Put $N(x):=\bigl(\log(109.2x)+M(1+\log2x)\bigr)/\log(1+\sqrt2)$. Since $N_s(x)-1\le N(x)$,
$N_0(x)-1\le N(x)$ and $E_B(k)(1+\log2x)\le N(x)$, each of the three inequalities of
Proposition~\ref{prop:X} implies
\[
x\log2<\log4+M\,N(x)\bigl(1+\log(N(x)x+2)\bigr).\tag{$*$}
\]
Let $x\ge M^3$ and $u:=\log x$, so that $u\ge145$ and $\log M\le u/3$. Then
$N(x)\le1.135\,M(u+4.7)$, and $1+\log(N(x)x+2)\le1.377\,u$ because $\log(u+4.7)\le0.035\,u$.
Hence the right-hand side of $(*)$ is at most $1.62\,M^2u^2\le1.62\,e^{2u/3}u^2<e^u\log2$, so $(*)$
fails. Thus $X(k)\le M^3<2.36\cdot10^{33}k^{14.4}$, and
\[
2.022\,E_k(X(k))\le G(k):=4.83\cdot10^{33}\,k^{14.4}(0.45k+2.02)\,\varphi^{-2k}.
\]
Since $14.4/k+0.45/(0.45k+2.02)<2\log\varphi$, the function $G$ is decreasing for $k\ge120$, and
$G(170)<10^{-3}$, which proves the claim for $k\ge170$. For $120\le k\le169$, a direct computation
of $X(k)$ (Appendix~A) shows that the largest value of $\log_{10}\bigl(2.022\,E_k(X(k))\bigr)$ is around
$-2.055$, attained at $k=120$.
\end{proof}


\begin{proposition}\label{prop:Kabs}
Equation \eqref{eq:main} has no solution satisfying \eqref{eq:remaining} with
$k\ge K_{\rm abs}:=1.5694\times10^{45}$.
\end{proposition}
\begin{proof}
Let $k\ge K_{\rm abs}$ and $x<X:=X(k)$; at $k=K_{\rm abs}$ we have $X<5.18\times10^{397}$. By
Lemma~\ref{lem:Esmall}, Lemma~\ref{lem:transfer} applies, and, comparing with Matveev's lower
bounds, the terms $\lambda_*$ and $\lambda_\zeta$ are never the minimum, because $2k\log\varphi$ is
much larger than the corresponding right-hand sides.

\emph{Case $d\ge1$.} If \eqref{eq:star} holds and the minimum in $\lambda_B$ is
$dx\log2-\log6.07$, then $dx<Y_L:=1.626\times10^{11}$. If moreover \eqref{eq:starr} holds, then
$\Phi_C$ gives $n<3.41\times10^{24}$; if \eqref{eq:starr} fails, then $r$ is bounded and $n=r+d<3.41\times10^{24}$
as well. Then $n\le k+1$ and $m\le(n-1)x+2<5.55\times10^{35}<k+1$, so that all the terms of
\eqref{eq:main} are odd-indexed Fibonacci numbers by \eqref{eq:kilic}, and \eqref{eq:main} becomes
\eqref{eq:fib}, which has no solution by Theorem~\ref{thm:fib}. Otherwise
$n<N_b:=3.86\times10^{12}$). Then $n\le k+1$, and
$\Phi_A$ gives $x\le dx<7.61\times10^{27}$, so that $m<2.93\times10^{40}<k+1$, and we conclude as before.

\emph{Case $d=0$.} Here $\Phi_{B0}$ gives $n<9.94\times10^{14}$, and then $\Phi_{A0}$ gives
$x<1.581\times10^{30}$, and the choice of $K_{\rm abs}$ gives $m\le(n-1)x+2\le K_{\rm abs}+1\le k+1$. This is the condition which determines
$K_{\rm abs}$. We conclude with Theorem~\ref{thm:fib}.

All the inequalities used are of the form $k>F(\log k)$ with $F$ a polynomial with positive
coefficients (because $X(k)<k^9$). Since they hold at $k=K_{\rm abs}$, they hold for all larger
$k$. We also checked them numerically for $k=K_{\rm abs}\cdot10^t$, $t=0,5,\ldots,180$.
\end{proof}

We now lower $K_{\rm abs}$ by reduction. Let $K$ be an upper bound for $k$ and consider
$120\le k\le K$, so that $x<X(K)$ and Lemma~\ref{lem:transfer} applies. One pass, with
$x\le M$ and $\log\alpha\ge2\log\varphi+\log(1-\varphi^{-30})$, consists of the following steps.
\begin{enumerate}[label=\rm(R\arabic*)]
\item ($d\ge1$.) Lemma~\ref{lem:leg} applied to $(x-1)\tau_0+e$, $\tau_0=\log\sqrt5/\log\varphi$,
with $A=2/\log\varphi$, gives $\lambda_B<W$. Hence either $dx<(W+\log6.07)/\log2$, or $n$ is
bounded, or $k$ is bounded (if $\lambda_*$ is the minimum), or \eqref{eq:star} fails.
\item ($d\ge1$, $dx<(W+\log6.07)/\log2$.) For each such $(d,x)$ and $x-2\le j\le2x-2$, $\Phi_C$
is explicit, and (b) bounds $n$ or $k$. If \eqref{eq:starr} fails, then $n=r+d$ is bounded.
\item ($d=0$.) Lemma~\ref{lem:dp} applied to $(x-1)\tau_0+e-\log2/\log\varphi$ with $A=2/\log\varphi$,
$B=e$, gives $\lambda_{B0}<W_0$, hence a bound for $n$ or for $k$.
\item For each $3\le n\le N$ (the resulting bound for $n$), Lemma~\ref{lem:dp} applied to $\Phi_A$ and $\Phi_{A0}$ gives a
bound for $x$, unless $k$ is smaller than an explicit bound.
\end{enumerate}
If $k$ exceeds all the bounds for $k$ met in (R1)--(R4) and $N$, we get $n\le N$ and
$x\le X'$. We repeat the pass with $M=X'$ until the bound stabilizes. Finally, if
$k+1\ge(N'-1)X'+2$ for the final bounds $N'$, $X'$, then $m\le k+1$ and $n\le k+1$, and
\eqref{eq:main} becomes \eqref{eq:fib}. This gives a new bound for $k$.

\begin{table}[ht]
\centering
\caption{Rounds of reduction for large $k$. Here $K$ is the bound for $k$ at the start of the
round; for the first pass we give $\lfloor\log_{10}X(K)\rfloor$, $W$ from (R1), the bound $N$ for
$n$ and the bound for $x$ from (R4). The last pass of every round gives
$W=8.28$, $N=15$ and $x\le43$, hence $(N-1)\cdot43+2=604$.}\label{tab:rounds}
\begin{tabular}{@{}rrrrrrr@{}}
\toprule
round & $K$ & $\lfloor\log_{10}X(K)\rfloor$ & $W$ & $N$ & $x\le$ & new $K$\\
\midrule
1 & $1.5694\times10^{45}$ & $397$ & $925.34$ & $1914$ & $1388$ & $2022$\\
2 & $2023$ & $56$ & $135.35$ & $278$ & $262$ & $604$\\
3 & $605$ & $52$ & $124.97$ & $256$ & $245$ & $604$\\
\bottomrule
\end{tabular}
\end{table}

\begin{proposition}\label{prop:rounds}
Equation \eqref{eq:main} has no solution satisfying \eqref{eq:remaining} with $k>K_0:=604$.
\end{proposition}
\begin{proof}
By Proposition~\ref{prop:Kabs} we may assume $k<K_{\rm abs}$.  Table~\ref{tab:rounds} shows the
rounds, each one starting from the bound given by the previous one, and the last one gives
again $k\le K_0$.
\end{proof}
\section{The cases $2\le k\le K_0$}\label{sec:smallk}

What remains is a finite computation, but not a small one: for each of the $603$ values of $k$ we
have a different $\alpha$, a different $g$, and a different bound $X(k)$, up to $10^{52}$. The
steps are those of Section~\ref{sec:fib}, in $\Q(\alpha)$ instead of $\Q(\sqrt5)$.

Let $2\le k\le K_0$. We start with $x<X(k)$ (Proposition~\ref{prop:X}) and repeat the following
pass, with $x\le M$, until the bound for $x$ stabilizes.
\begin{enumerate}[label=\rm(S\arabic*)]
\item ($d\ge1$.) With $\tau:=-\log g/\log\alpha$, Lemma~\ref{lem:Lambda}(b) gives
$|(x-1)\tau-j|<(6/\log\alpha)e^{-w}$ with $w=\min\{dx\log2,\ n\log\alpha-\log x\}$, and
Lemma~\ref{lem:leg} with $u=x-1\le M$ gives $w<W$. Hence either $dx<W/\log2$, or
$n<(W+\log M)/\log\alpha$, or \eqref{eq:star} fails and $n<1+\log(36.4M)/\log\alpha$.
\item ($d\ge1$, $dx<W/\log2$.) For each such $(d,x)$ and $x-2\le j\le2x-2$ we computed
$\Lambda_C$; when it is non-zero, Lemma~\ref{lem:Lambda}(c) gives $n<\log(2.7x/|\Lambda_C|)/
\log\alpha$, and if \eqref{eq:starr} fails, then $n=r+d<1+\log(36.4x)/\log\alpha+d$. The only
zero is $(k,d,x,j)=(2,1,2,1)$ (Lemma~\ref{lem:zero}), which is excluded by
\eqref{eq:remaining}; the smallest non-zero value is $|\Lambda_C|=8.9\times10^{-6}$.
\item ($d=0$.) Lemma~\ref{lem:dp} applied to $(x-1)\tau-j-\log2/\log\alpha$, with $u=x-1$,
$A=3.8M/\log\alpha$, $B=\alpha$ and $w=n$, bounds $n$.
\item For each $3\le n\le N$ (the bound from (S1)--(S3)), Lemma~\ref{lem:dp} applied to
$\Lambda_A$ (with $\tau=\log\alpha/\log P_n$, $\mu=\log g/\log P_n$, $u=m\le(n-1)M+2$, $A=3/\log P_n$,
$B=2$, $w=dx$) and to $\Lambda_{A0}$ ($\mu=\log(g/2)/\log P_n$, $A=0.5/\log P_n$, $B=5$, $w=x$)
gives a new bound for $x$. For large $n$, $\log g$ is very close to $\log P_n-n\log\alpha$ and
Lemma~\ref{lem:dp} fails for $\Lambda_A$; we then use Lemma~\ref{lem:leg} with
$\tau_n=\log P_n/\log\alpha$, $u=x-1$ and $\delta_n=|\log(P_n/(g\alpha^n))|/\log\alpha$, exactly as
in Section~\ref{sec:fib}.
\end{enumerate}

\begin{table}[ht]
\centering
\caption{Maximum over $2\le k\le K_0$ of the quantities met in the passes (for most $k$ the bound
stabilizes after three or four passes; a fifth pass was needed for $25$ values of $k$, and a sixth one for one value).}\label{tab:smallk}
\begin{tabular}{@{}rrrrrr@{}}
\toprule
pass & $M$ & $W$ in (S1) & $dx<W/\log 2$ & $N$ & new bound for $x$\\
\midrule
1 & $1.44\times10^{52}$ & $125.86$ & $181$ & $255$ & $256$\\
2 & $257$ & $13.94$ & $20$ & $20$ & $53$\\
3 & $54$ & $10.28$ & $14$ & $17$ & $51$\\
4 & $48$ & $10.18$ & $14$ & $17$ & $51$\\
5 & $46$ & $10.14$ & $14$ & $16$ & $50$\\
\bottomrule
\end{tabular}
\end{table}

\begin{proposition}\label{prop:smallk}
Let $2\le k\le K_0$. Every solution of \eqref{eq:main} satisfying \eqref{eq:remaining} has
$n\le17$ and $x\le51$, and there is no such solution.
\end{proposition}
\begin{proof}
The bounds follow from Table~\ref{tab:smallk}. For each $k$, a search in exact integer arithmetic
over $3\le n\le N_k$, $1\le r\le n$, $2\le x\le X_k$, where $N_k$ and $X_k$ are the final bounds for
this $k$, comparing $P_n^x+P_r^x$ with the list of the $P_m$, $m\le(N_k-1)X_k+3$, finds only, for
$k=2$, the solutions $(n,n-1,2n-1,2)$, which are excluded by \eqref{eq:remaining}.
\end{proof}

\section{Proof of Theorem~\ref{thm:main}}\label{sec:proof}

Let $(k,n,r,m,x)$ be a solution of \eqref{eq:main} with $n\ge r$. If $x=1$, or $n\le2$, or
$(k,x,d)=(2,2,1)$, Lemmas~\ref{lem:x1}, \ref{lem:n2} and \ref{lem:pell} give the solutions listed
in Theorem~\ref{thm:main}. Otherwise \eqref{eq:remaining} holds. By
Proposition~\ref{prop:rounds} we have $k\le K_0$, and Proposition~\ref{prop:smallk} shows that
there is no solution. This proves Theorem~\ref{thm:main}.\qed

\appendix
\section{Computational details}\label{app:comp}

All computations were done in Python~3, with the \texttt{decimal} module for real numbers and exact
integer arithmetic for the sequences. The precision was chosen in each step as twice the number of
digits of the bound $M$ plus $80$ to $140$ digits, so that the products $\tau q$ and $\mu q$ in
Lemma~\ref{lem:dp}, with $q$ up to $6M\cdot10^{20}$, are correct. The dominant root $\alpha$ was
computed by Newton's method on $x^2-3x+1+x^{1-k}=0$, which is well conditioned for large $k$. The
discriminants $\Delta_k$ of Lemma~\ref{lem:disc} were also computed exactly for $2\le k\le1500$,
from the formula in the proof (reduction of $\Psi_k$ modulo $q$ in exact rational arithmetic), and
compared with a floating-point evaluation of $\prod_i\Psi_k'(\alpha_i)$ for $k\le8$.

As a consistency check, independent of the proof, a direct search over $2\le k\le80$,
$1\le r\le n\le45$ and $1\le x\le25$, and over $2\le k\le12$, $1\le r\le n\le60$ and
$1\le x\le30$, found only the solutions of Theorem~\ref{thm:main}.

The scripts are: \texttt{fibpart.py} (Section~\ref{sec:fib}), \texttt{perk.py}
(Proposition~\ref{prop:X} and Section~\ref{sec:smallk}), \texttt{largek.py}
(Section~\ref{sec:largek}), \texttt{disc.py} (Lemma~\ref{lem:disc}) and \texttt{verify\_all.py},
which repeats all the steps with assertions. On one core, Section~\ref{sec:fib} takes a few seconds, the rounds of Section~\ref{sec:largek} about $40$ minutes (the first pass works with about $950$ digits), Section~\ref{sec:smallk} about $30$ minutes, and the complete verification script about $70$ minutes.

\section*{Declarations}

\paragraph{\textbf{Data availability}} The Python scripts used in this paper are available from
the second author on request.
interest.

\paragraph{\textbf{Use of AI tools}} Claude (Anthropic) was used to search the literature,  to write and run the Python
code described in Appendix~\ref{app:comp}, and for typeediting. The authors take full responsibility for the content.

\end{document}